\documentclass[10pt,reqno]{amsart}
\usepackage{graphicx}
\usepackage{indentfirst,csquotes}

\usepackage{amssymb,amsthm,amsmath}
\usepackage{xcolor,paralist,hyperref,titlesec,fancyhdr,etoolbox}

\usepackage{titlesec}

\hypersetup{ colorlinks=true, linkcolor=black, filecolor=black, urlcolor=black }

\RequirePackage{fix-cm}

\usepackage[english]{babel}

\usepackage{graphicx}
\usepackage{tabularx}
\usepackage{color}
\usepackage{mathtools,amsfonts,amssymb,amsmath,bm}
\usepackage{url}            
\usepackage{soul}
\renewcommand{\d}{{\rm d}}

\begin{document}
\title[Dynamics of a mathematical model of virus spreading]{Dynamics of a mathematical model of virus spreading incorporating the effect of a vaccine} 
\pagestyle{plain}
\author{Ayt\"ul G\"ok\c{c}e{$^{1}$} \and Burcu G\"urb\"uz{$^{*, 2}$} \and Alan D. Rendall{$^{3}$} }
\date{\today}
\address{A. G\"ok\c{c}e{$^{1}$} Ordu University, Faculty of Science and Letters, Department of Mathematics, 52200, Ordu, Turkey}
\email{aytulgokce@odu.edu.tr}

\address{{B. G\"urb\"uz$^{*, 2}$} Institut f\"ur Mathematik, Johannes Gutenberg-Universit\"at, Staudingerweg 9, 55099, Mainz, Germany}
\email{burcu.gurbuz@uni-mainz.de}

\address{A. D. Rendall{$^{3}$} Institut f\"ur Mathematik, Johannes Gutenberg-Universit\"at, Staudingerweg 9, 55099, Mainz, Germany}
\email{rendall@uni-mainz.de}

\let\thefootnote\relax
\footnotetext{MSC2020: 00A71, 34D20, 37M05, 37N25, 92D30.} 

\begin{abstract}
The COVID-19 pandemic led to widespread interest in epidemiological models.
In this context the role of vaccination in influencing the spreading of the
disease is of particular interest. There has also been a lot of debate on the
role of non-pharmaceutical interventions such as the disinfection of surfaces.
We investigate a mathematical model for the spread of a disease which includes
both imperfect vaccination and infection due to virus in the environment. The
latter is studied with the help of two phenomenological models for the force
of infection. In one of these models we find that backward bifurcations take
place so that for some parameter values an endemic steady state exists
although the basic reproduction ratio $R_0$ is less than one. We also prove
that in that case there can exist more than one endemic steady state. In the
other model all generic transcritical bifurcations are forward bifurcations so
that these effects cannot occur. Thus we see that the occurrence of backward
bifurcations, which can be important for disease control strategies, is
dependent on the details of the function describing the force of infection. By
means of simulations the predictions of this model are compared with data for
COVID-19 from Turkey. A sensitivity analysis is also carried out.
\end{abstract} 
\maketitle

\textbf{1. Introduction}\\

\noindent Since the beginning of epidemiology mathematical models have played a central
role. This can be seen in the groundbreaking work of Ronald Ross and Hilda
Hudson on the eradication of malaria \cite{ross16}, \cite{ross17a},
\cite{ross17b}. In that work the authors identified a threshold for the
persistence of a disease which can be seen as the ancestor of the basic
reproduction ratio $R_0$ which is so important in epidemiology today. The
COVID-19 epidemic caused a surge of work where epidemiological models were
defined, simulated and subjected to rigorous mathematical analysis. Due to
the urgency of the situation this development took place in a rather
disorganized way. Now it is time to consolidate and extend the things learned
at that time so as to be prepared as well as possible for future epidemics.

\noindent In this paper we study a model for the spread of an infectious disease
in a human population which includes an imperfect vaccination and takes into account infections due to virus particles in the environment. In
particular we are thinking of fomites, objects in the environment which are contaminated with virus and which are not humans or animals. Here we may think of the contamination by hands touching doorknobs \cite{zhao12} or infections spreading in hospitals \cite{sulyok21}. The question of how important this route of infection is for COVID-19 has been a subject of much discussion. The consensus appears to be that it is of secondary importance but this may be different for other diseases \cite{anderson20,park20}.

\noindent In Section 2 
we define the model which is of central
interest in this paper and establish some basic properties of its solutions.
The model contains a response function which describes how the concentration
of virus in the environment affects the rate of infection by this route. This
function depends on an integer $n\ge 1$. The motivation for the choice of this
function is also discussed. In Section 2.1 
it is shown that this
model has a unique disease-free steady state (Lemma 3.1) and the stability of
that state is determined using the next generation matrix. The basic
reproduction ratio $R_0$ is computed for this model and it is shown to be
a decreasing function of the vaccination rate. The existence of backward
bifurcations is analysed using the method of van den Driessche and Watmough
\cite{van2002reproduction}.  It is proved that generic backward bifurcations
occur in the case $n=2$ but not otherwise (Theorem 3.1). In particular they do
not occur in the case $n=1$, where the
function describing the force of infection is one which had previously been
considered in the literature \cite{bulut2021modelling}. Note that an existing
model for imperfect vaccination \cite{gumel2006sveir} also exhibits no
backward bifurcations. A backward bifurcation
is often accompanied by the presence of more than one positive steady state
for given values of the parameters (cf. \cite{martcheva19}).
In many cases of backward bifurcations
simulations show not only that for a certain choice of parameters positive
steady states exist for $R_0<1$ but also that two positive steady
states can occur. It is proved in Section 4 
that there are parameters for which our model with $n=2$ exhibits the latter behaviour.
In Section 5 
it is shown that solutions of the model can be fitted to COVID-19 data from Turkey. Section 6 
carries out a sensitivity analysis of the model.

\textbf{2. The model}

The model considered in what follows is a generalization of one introduced in
\cite{gumel2006sveir} to study the effects of vaccination against SARS and is
given by the following equations:
\begin{align}
    \frac{d S}{d t} &= \Lambda-\beta SI - \sigma S+(1-\lambda)t'V-\alpha_1 S g(C,\kappa)-\mu S, \label{Eq:ModelExp1}\\
    \frac{d E}{d t}&= \beta SI +\epsilon I V - \xi E-\mu E +\alpha_1 S g(C,\kappa)+ \alpha_2 V g(C,\kappa), \label{Eq:ModelExp2}\\
    \frac{d I}{d t} &= \xi E -\delta I - d I -\mu I, \label{Eq:ModelExp3}\\
    \frac{d V}{d t} &= \sigma S-\epsilon I V - (t'+\mu) V-\alpha_2 V g(C,\kappa), \label{Eq:ModelExp4}\\
    \frac{d R}{d t} &= \delta I - \mu R+\lambda t' V, \label{Eq:ModelExp5}\\
    \frac{d C}{d t} &= \varphi I -\omega C, \label{Eq:ModelExp6}
\end{align}
where \[g(C,\kappa)=\frac{C^n}{C^n+\kappa}.\]
The meaning of the parameters in this model is described in Table 1. The model
of \cite{gumel2006sveir} was called an SVEIR model after the names of its five
unknowns $S$, $V$, $E$, $I$ and $R$. These are the numbers of susceptible,
vaccinated, exposed, infectious and recovered individuals, respectively. We
augment this by an additional variable $C$ representing the concentration of
the virus in the environment. In both models the vaccination is imperfect but
the imperfection is of a different kind. Correspondingly the class $V$ has a
different interpretation in the two models. In \cite{gumel2006sveir} the class
$V$ consists of individuals who have been vaccinated at some time. The effect
of the vaccination is to lower the rate at which they get infected compared to
unvaccinated individuals. In the model
(\ref{Eq:ModelExp1})-(\ref{Eq:ModelExp6}) the class $V$ consists of
individuals who have received a vaccination but where the vaccination has not
yet had time to become fully effective. After that time either the
vaccination provides complete protection or it has not been effective and the
individual returns to the susceptible class. For biological reasons the
inequalities $\epsilon\le\beta$ and $\alpha_2\le\alpha_1$ are assumed, which
means that vaccinated individuals are no more likely to be infected than
unvaccinated individuals, either by infected individuals or by contact with their surroundings.

\noindent Mathematically the model of \cite{gumel2006sveir}, up to a different notation,
can be obtained from our model by setting $t'=\alpha_1=\alpha_2=0$ and
discarding the equation for $S$. This is possible since when the parameters
just listed are zero the equations for the first five variables to not depend
on $S$. In the model (\ref{Eq:ModelExp1})-(\ref{Eq:ModelExp6}) the
imperfection of the vaccination is expressed as follows. Individuals leave the
vaccinated state at rate $t'$, the vaccination having been successful with
probability $\lambda$. This way of modelling an imperfect vaccination was
previously used in \cite{angeli2022modeling}. The other additional effect
taken account of in (\ref{Eq:ModelExp1})-(\ref{Eq:ModelExp6}) is related to
infection by virus in the environment. It is expressed by the terms containing
the factors $\alpha_1$ and $\alpha_2$ relating to the unvaccinated and
vaccinated individuals, respectively. This type of effect was included in a
model of \cite{bulut2021modelling}. In that paper the function $g$ written
above with $n=1$ was used as a phenomenological description of the rate of infection in this process.

\noindent It is worth taking some time to discuss the status of this type of
phenomenological description. It is used in defining response functions in
various parts of biology. In biochemistry the case $n=1$ is called a
Michaelis-Menten function while the case $n\ge 2$ is called a Hill function. In
predator-prey models in ecology the case $n=1$ is called Holling type II while
the case $n\ge 2$ is called Holling type III. Holling type I denotes a linear
response function, usually with a cut-off. A general discussion of response
functions in epidemiology is given in Chapter 10 of the book of Diekmann and
Heesterbeek \cite{diekmann00}, whereby the authors make clear from the
beginning that they do not claim to give a definitive answer to the questions they are raising.

\noindent Suppose that there is a source of infection with intensity $Z$ and a population
$S$ of susceptibles. Let $F(Z,S)$ be the rate of infection. In principle this
could be any function. Let us suppose that $F$ depends linearly on $S$ but
initially allow its dependence on $Z$ to be arbitrary. Thus $F(Z,S)=sf(Z)$ for
some function $f$. What properties should the function $f$ have? It should
be positive for $Z$ positive and zero for $Z=0$. It should be non-decreasing.
It is reasonable to assume that it is bounded. The simplest type of
function satisfying these requirements is one of the form $f(Z)=\frac{aZ}{Z+b}$.
Another situation in which a response function is of relevance is the predation
rate in a predator-prey model. There the analogue of $I$ is the density of
predators while the analogue of $Z$ is the density of prey. In that case a
function of the form just considered is called Holling type II. In that context
Holling type I is a function which is linear up to a threshold value and then
constant. Holling type III corresponds to $Z$ being replaced by $Z^p$. This
argument for introducing a function $f$ of this form is purely
phenomenological. Holling had a mechanistic argument to motivate his
type II function. There have also been attempts to motivate the type III
function by mechanistic considerations (cf. \cite{dawes13}, \cite{holling59}).
We are not aware that this has been done in epidemiology. The function
corresponding to Holling type II was introduced to epidemiological models by
Dietz \cite{dietz82}, without a mechanistic background.
Holling's mechanistic approach does not apply to epidemiological models.
Diekmann and Heesterbeek \cite{diekmann00} discuss mechanistic approaches to
the Holling type II function in epidemiology. In fact in a model case they
derive something which is not a rational function. We have not found a paper
where the Holling type III function is used in epidemiology.

\begin{table}[h!]
\caption{ Parameters used in the model. \label{Tab:11}}
\begin{center}
\begin{tabular}{ c l }
\hline \hline
\small\textbf{parameters}  &  \small\textbf{biological meaning}\\
\hline \hline
$\Lambda$  & recruitment rate\\
$\beta$ & effective contact rate with $\beta S$ new  susceptible  individuals per unit time\\
$\alpha_{1}$ & transmission ratio of the virus from the environment  to susceptible individuals\\& that enter the exposed class\\
$\alpha_{2}$ & transmission ratio of the virus from the environment  to vaccinated individuals\\& (not fully immunised) that may enter the exposed class\\
$\epsilon$ &  the rate at which a vaccinated individual (not fully immunised) becomes exposed \\& after being in contact with an infected individual\\
$\delta$ &  the rate of recovered  individuals\\
$\mu$  & natural mortality rate\\
$\xi$ & rate of development of clinical symptoms\\
$ d$  & disease induced fatality rate\\
$\sigma$ &  vaccination rate of susceptible individuals (the first shot)\\
$\varphi$ & the virus exposure rate \\
$\lambda$ & the efficiency of the  vaccine \\
$t'^{-1}$ & the mean amount of time that is spent in the vaccinated class before developing \\& an immune response and moving to the recovery class.\\
$\omega$ & the rate of  decay in the virus density\\
\hline
\end{tabular}
\end{center}
\end{table}

\noindent The right hand sides of equations \eqref{Eq:ModelExp1}-\eqref{Eq:ModelExp6}
are smooth and hence for any initial values at a given time they have a unique
solution on some time interval. Because of the interpretation of the unknowns
we are interested in solutions which are non-negative at all times. This is
true provided the initial data are non-negative (cf. \cite{rendall12},
Lemma 1). Let $N(t) = S(t)+E(t)+I(t)+V(t)+R(t)$. Then
\begin{equation}
    \frac{\d N}{\d t} = \Lambda -\mu N - d I\le \Lambda -\mu N.
\end{equation}

\noindent This implies that $N$ remains bounded on any finite time interval. Hence on
such an interval all variables other than $C$ are bounded. It then follows
from \eqref{Eq:ModelExp6} that $C$ is also bounded. As a consequence the
solutions exist globally in the future. Using the differential inequality for
$N$ again shows that
\begin{equation}
   \limsup\limits_{t \to \infty} N(t) \leq \frac{\Lambda}{\mu}.
\end{equation}
It then follows from \eqref{Eq:ModelExp6} that
$\limsup_{t\to\infty} C(t)\le\frac{\Lambda\varphi}{\mu\omega}$.

\textbf{2.1. The disease-free steady state}

Consider a boundary steady state $(S^*,E^*,I^*,V^*,R^*,C^*)$ of the system
(\ref{Eq:ModelExp1})-(\ref{Eq:ModelExp6}), i.e. a time-independent solution
for which at least one of the unknowns is zero.

\noindent
{\bf Lemma 2.1} The model (\ref{Eq:ModelExp1})-(\ref{Eq:ModelExp6}) has a unique
boundary steady state.

\noindent
{\bf Proof} Let $(S^*,E^*,I^*,V^*,R^*,C^*)$ be a boundary steady state. If
$S^*=0$ then (\ref{Eq:ModelExp1}) gives a contradiction and so $S^*\ne 0$. If
$V^*=0$ then (\ref{Eq:ModelExp4}) implies that $S^*=0$. Hence in fact
$V^*\ne 0$. If $R^*=0$ then (\ref{Eq:ModelExp6}) implies that $V^*=0$. Hence
in fact $R^*\ne 0$. It follows from the other three equations that for a
steady state the equations $E^*=0$, $I^*=0$ and $C^*=0$ are all equivalent to
each other. Hence any boundary steady state is of the form
$(S^*,0,0,V^*,R^*,0)$.
In this case the steady state equations are equivalent to the system
\begin{align}
 \Lambda-\sigma S^*+(1-\lambda)t'V^*-\mu S^*=0, \nonumber \\
 \sigma S^*-(t'+\mu)V^*=0, \label{Eq:DFExp2}\\
 \mu R^*+\lambda t'V^*=0,\nonumber
\end{align}
\noindent and these can be solved to give
\begin{align}
S^*&=\frac{\Lambda (t'+\mu)}{\mu(\sigma+t'+\mu)+\lambda t' \sigma }, \label{Equ:Ss}\\
V^*&=\frac{\Lambda \sigma }{\mu(\sigma+t'+\mu)+\lambda t'\sigma},  \label{Equ:Vs}\\
R^*&=\frac{\lambda t' \Lambda \sigma}{\mu(\mu(\sigma+t'+\mu)+\lambda t'\sigma)}.  \label{Equ:Rs}
\end{align}
\noindent Thus there exists a unique boundary steady state of this model. $\blacksquare$

\noindent All the variables corresponding to the presence of infection are zero in this
state and so we call it the disease-free steady state.

\begin{figure}[ht!]
\centering
\includegraphics[width=0.6\textwidth]{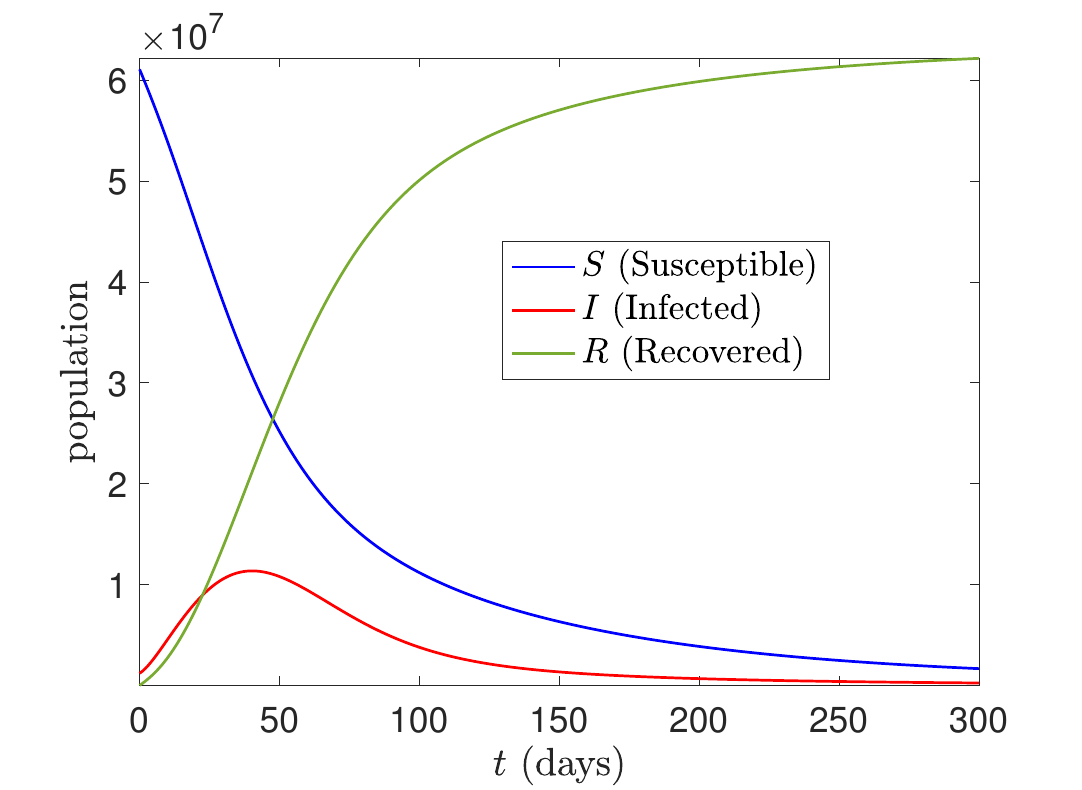} \\
\caption{Simulation result of the model \eqref{Eq:ModelExp1}-\eqref{Eq:ModelExp6}
with initial data and parameters given in the text. }
\label{Fig:SIRV}
\end{figure}

\noindent To illustrate how a solution of this model corresponding to an epidemic might
look we show the results of a simulation with biologically motivated
parameters. Table 2 lists the references which were used either as direct
sources or guidelines for the choice of the parameters.
Figure \ref{Fig:SIRV} demonstrates the dynamics of populations during an epidemic for $300$ days. Susceptible ($S$), Infected ($I$) and Recovered ($R$) populations respectively associated with blue, red and green lines. The initial conditions are chosen as $S_0=61098000, V_0=18500000, E_0=2200000, I_0= 1200000$ and $R_0=2000$ and parameters are given in Table \ref{tab:Hypo}. It should be noted that these initial conditions and parameters are only selected for illustrative purposes and may not be epidemiologically realistic.

\begin{table}[]
    \centering
       \caption{Biologically meaningful parameters used in Fig. \ref{Fig:SIRV}.}
    \begin{tabular}{c l l l}
       \hline
       Parameters  &  Value ($n=1$)  & Unit  & Source\\
       \hline
       $\Lambda$ & $3032$ & day$^{-1}$ & assumed based on  \cite{bajiya2020mathematical,garba2013cross,gumel2006sveir} \\
        $\beta$ & $0.15\times 10^{-8}$   & day$^{-1}$ & assumed based on \cite{gumel2006sveir} \\
        $\mu$ & $3.653 \times 10^{-5} $   & day$^{-1}$ & assumed based
 on \cite{bulut2021modelling,gumel2006sveir} \\
        $\epsilon$ & $0.15\times 10^{-8}$ &  day$^{-1}$ & assumed based on \cite{angeli2022modeling,gumel2006sveir}  \\
        $t'$ & $1/120 $   & day$^{-1}$ & assumed based on \cite{angeli2022modeling} \\
        $\lambda$ & $0.8$  & day$^{-1}$ & assumed based on \cite{angeli2022modeling,baden2021efficacy,polack2020safety}  \\
        $d$     &  $0.02$    & day$^{-1}$ & assumed based on \cite{bulut2021modelling,ritchie2020mortality} \\
        $\alpha_1$ &  $0.01$ & day$^{-1}$ & assumed  \\
        $\alpha_2$ & $0.01$  & day$^{-1}$ & assumed   \\
        $\omega$ & $4$ & day$^{-1}$ &  \cite{bulut2021modelling}  \\
        {$\varphi$} & $2$ & day$^{-1}$ &  \cite{bulut2021modelling} \\
        $\kappa$ & $20000$ & copies $/$day & assumed based on \cite{bulut2021modelling} \\
        $\xi$ & {$0.125$} & day$^{-1}$ & \cite{gumel2006sveir} \\
        $\delta$ & $0.06$ & day$^{-1}$ &  assumed based on \cite{angeli2022modeling,gumel2006sveir}\\
        $\sigma$ & $0.01$  & day$^{-1}$ & \cite{angeli2022modeling} \\
    \end{tabular}
    \label{tab:Hypo}
\end{table}

\noindent Here a dramatic increase  can be seen in the number of infected individuals until day $45$, then a gentle  decline appears for the population of infected individuals. As is observed from the graph the number of susceptible individuals slowly decreases to $16000000$ and  while the number of recovered individuals rises above $62000000$ at day $290$. The total population is taken as $83$ million.

\textbf{3. Stability of the disease-free steady state}

\noindent Linearisation around the disease-free steady state leads to the Jacobian matrix
\begin{equation}
    \label{Eq:JacobDef}
    \mathcal{J}=[J_{ij}]_{6\times6}\biggr\rvert_{E^*}, \ \mbox{for} \ i,j=1,2,...,6,
\end{equation}
where $E^*=(S^*,0,0,V^*,R^*,0)$ and
\begin{align*}
   J_{11}&=-\sigma-\mu, \  J_{12}=0, \  J_{13}=-\beta S, \  J_{14}=(1-\lambda)t',\  J_{15}=0,\  J_{16}=-{\alpha_1 S}\delta_{1n}/{\kappa},\\
   J_{21}&=0, \  J_{22}=-\xi-\mu, \  J_{23}=\epsilon V+\beta S, \  J_{24}=0,\  J_{25}=0,\  J_{26}=(\alpha_1 S + \alpha_2 V)\delta_{1n}/\kappa,\\
   J_{31}&=0, \  J_{32}=\xi, \  J_{33}=-(\delta +d+\mu), \  J_{34}=0,\  J_{35}=0,\  J_{36}=0,\\
   J_{41}&=\sigma, \  J_{42}=0, \  J_{43}=-\epsilon V, \  J_{44}=-(t'+\mu),\  J_{45}=0,\  J_{46}=-\alpha_2 V\delta_{1n}/\kappa,\\
   J_{51}&=0, \  J_{52}=0, \  J_{53}=\delta, \  J_{54}=\lambda t',\  J_{55}=-\mu, J_{56}=0,\\
   J_{61}&=0, \  J_{62}=0, \  J_{63}=\varphi, \  J_{64}=0,\  J_{65}=0,\  J_{66}=-\omega.\\
\end{align*}
Here $\delta_{1n}$ is a Kronecker delta. Then we have the matrix
\begin{align*}
\mathcal{J}&=
  \begin{bmatrix}
   J_{11} & 0 & J_{13} & J_{14} & 0 & J_{16} \\
   0 & J_{22} & J_{23} & 0 & 0 & J_{26} \\
   0 & J_{32} & J_{33} & 0 & 0 & 0 \\
   J_{41} & 0 & J_{43} & J_{44} & 0 & J_{46} \\
   0 & 0 & J_{53} & J_{54} & J_{55} & 0 \\
   0 & 0 & J_{63} & 0 & 0 & J_{66} \\
\end{bmatrix}
\end{align*}
\noindent It is clear that one of the eigenvalues of the Jacobian is $J_{55}$. Moreover,
removing the fifth row and column and interchanging the second row and column
with the fourth leads to a matrix with block diagonal structure. Thus two
further eigenvalues of the Jacobian can be obtained as the eigenvalues of
the matrix
\begin{align*}
  \overline{\mathcal{J}_1}=
    \begin{bmatrix}
   J_{11} & J_{14} \\
   J_{41} & J_{44} \\
\end{bmatrix}.
\end{align*}
\noindent Now
\begin{align*}
{\rm tr\ } \overline{\mathcal{J}_1}  &= J_{11}+J_{44}=-(\sigma+\mu)-(t'+\mu)<0\\
    \det \overline{\mathcal{J}_1}&= (\sigma +\mu)(t'+\mu)-\sigma (1-\lambda)t'\\
        &= \mu(\sigma +\mu + t')+\sigma \lambda t' >0.
\end{align*}
\noindent Thus the eigenvalues of $\overline{\mathcal{J}_1}$ have negative real parts. The
remaining three eigenvalues of the Jacobian are the eigenvalues of the matrix
\begin{align*}
  \overline{\mathcal{J}_2}=
    \begin{bmatrix}
   J_{22} & J_{23} & J_{26}  \\
   J_{32} & J_{33} & 0 \\
   0 & J_{63} & J_{66}
\end{bmatrix} ,
\end{align*}
leading to characteristic polynomial
\begin{align*}
(\ell-J_{22})(\ell-J_{33})(\ell-J_{66})-J_{23} J_{32} (\ell-J_{66})- J_{26} J_{32} J_{63}=0,
\end{align*}
which can be rewritten as
\begin{align}
\ell^3+\mathcal{A}_1 \ell^2 + \mathcal{A}_2 \ell + \mathcal{A}_3 =0, \label{Equ:chace}
\end{align}
where
\begin{align}
\mathcal{A}_1 &=- J_{22}- J_{33}-J_{66}, \nonumber\\
\mathcal{A}_2 &= J_{22} J_{33}-J_{23}J_{32}+J_{66}(J_{22}+J_{33}),\nonumber\\
\mathcal{A}_3 &=- J_{22} J_{33} J_{66} \left(1-\frac{J_{32} \left(J_{23} J_{66} - J_{63} J_{26}\right)}{J_{22} J_{33} J_{66} }\right).\label{Equ:R01}
\end{align}
\noindent The Routh-Hurwitz criterion says that all roots of the characteristic equation
\eqref{Equ:chace} have negative real parts if and only if
$\mathcal{A}_1>0$, $\mathcal{A}_1 \mathcal{A}_2> \mathcal{A}_3$ and
$\mathcal{A}_3>0$ and if these conditions hold the disease-free steady state
is asymptotically stable. It is clear that the first condition holds but it is
not so easy to see when the second and third conditions hold. It will later
be proved indirectly using the next generation matrix that they hold in this
model for all values of the parameters.

\textbf{3.1. The next generation matrix}

In this section, following the ideas presented in \cite{van2002reproduction},
the basic reproduction ratio for (\ref{Eq:ModelExp1})-(\ref{Eq:ModelExp6})
is derived using the next generation matrix method.
We use the notation of \cite{van2002reproduction}. To apply this method we
must choose which of the unknowns represent groups of infected individuals and
which terms in the equations represent new infections. In fact we choose $E$,
$I$ and $C$ to be the infected variables and the terms which are non-negative
and non-linear in the unknowns to represent new infections. The conditions
(A1)-(A5) of \cite{van2002reproduction} are satisfied. Most of these are
rather obvious for this model. The only exception is (A5) which holds
because the quantities corresponding to $J_{23}$ and $J_{26}$  are zero in the
case that new infections are turned off. The matrix ${\mathcal {F}}$
associated with new infections and the matrix ${\mathcal {V}}$ containing
the remaining expressions  are given by
\begin{align*}
    {\mathcal {F}}&=
        \begin{bmatrix}
        0\\
        \beta SI +\epsilon VI +(\alpha_1 S+\alpha_2 V)g(C,\kappa) \\
        0\\
        0\\
        0\\
        0\\
        \end{bmatrix},
 \end{align*}
 and
 \begin{align*}
            -{\mathcal{V}}&=
                \begin{bmatrix}
        0\\
         \xi E+\mu E \\
        -\xi E+\delta I +dI+\mu I\\
        0\\
        0\\
        -\varphi I +\omega C\\
        \end{bmatrix},
        \end{align*}
${\mathcal{V}}_+$ and ${\\mathcal {V}}_-$ are the positive and negative parts of $\mathcal {V}$,
respectively. Hence the matrices $F$ and $V$ of \cite{van2002reproduction} are
given by
        \begin{align*}
    F\bigr\rvert_{E^*}&=
        \begin{bmatrix}
        0 & \beta S^*+\epsilon V^* &  \delta_{1n}\left(\alpha_1 S^* + \alpha_2 V^*\right)/\kappa \\
        0 & 0 & 0\\
        0 & 0 & 0\\
        \end{bmatrix}, \
    V\bigr\rvert_{E^*}=
        \begin{bmatrix}
        \xi+\mu & 0                 & 0 \\
        -\xi    & \delta+d+\mu      & 0\\
        0       & -\varphi          & \omega\\
        \end{bmatrix},
\end{align*}
The reproduction ratio $R_0$ is defined (cf. \cite{bulut2021modelling},
\cite{van2002reproduction}) to be the spectral radius of the matrix given by
\begin{align*}
    FV^{-1}=
        \begin{bmatrix}
         FV^{-1}_{11}      & FV^{-1}_{12} & FV^{-1}_{13} \\
        0       & 0                     & 0\\
        0       & 0                     & 0\\
        \end{bmatrix},
\end{align*}
where
\begin{align*}
 FV^{-1}_{11} &=  \frac{\xi}{(\xi+\mu)(\mu+\delta+d)}\left[\left( \beta + \frac{\delta_{1n}\alpha_1 \varphi}{\omega \kappa}\right)S^*+ \left( \epsilon + \frac{\delta_{1n}\alpha_2 \varphi}{\omega \kappa}\right)V^*\right],\\
 FV^{-1}_{12} &=   \frac{1}{(\mu+\delta+d)}\left[\left( \beta + \frac{\delta_{1n}\alpha_1 \varphi}{\omega \kappa}\right)S^*+ \left( \epsilon + \frac{\delta_{1n}\alpha_2 \varphi}{\omega \kappa}\right)V^*\right],\\
 FV^{-1}_{13} &=  \frac{\delta_{1n}}{\omega \kappa} \left(\alpha_1 S^* +\alpha_2 V^*\right)
\end{align*}
The characteristic equation of this matrix is given by
\begin{align*}
    \det(FV^{-1}-\Sigma I)=0.
\end{align*}
Its roots are the eigenvalues:
\begin{align}
    \label{Eq:Eigenvalues}
    \Sigma_1&=\frac{\xi}{(\xi+\mu)(\mu+\delta+d)}\left[\left( \beta + \frac{\delta_{1n}\alpha_1 \varphi}{\omega \kappa}\right)S^*+ \left( \epsilon + \frac{\delta_{1n}\alpha_2 \varphi}{\omega \kappa}\right)V^*\right], \\
    \Sigma_2&=0,\\
    \Sigma_3&=0.
\end{align}
Thus the basic reproduction ratio, which is associated with the dominant
eigenvalue $\Sigma_1$, is
\begin{align}
    R_0&=\frac{\xi S^*}{(\xi +\mu)(\mu+\delta+d)}\left[\left( \beta + \frac{\delta_{1n}\alpha_1 \varphi}{\omega \kappa}\right)+ \left( \epsilon + \frac{\delta_{1n}\alpha_2 \varphi}{\omega \kappa}\right) \frac{\sigma}{t'+\mu}\right],\label{Eq:R03}\\
    &= \frac{\xi}{\omega(\xi +\mu)(\mu+\delta+d)} \frac{\Lambda (t'+\mu)}{\mu (\sigma+t'+\mu)+\lambda t' \sigma} \left[\omega \beta +\frac{\delta_{1n}\alpha_1 \varphi}{\kappa}+\frac{\sigma}{t'+\mu}\left(\epsilon \omega+\frac{\delta_{1n}\alpha_2 \varphi}{\kappa}\right)\right]. \nonumber
\end{align}
It follows from Theorem 2 of \cite{van2002reproduction} that the disease-free
steady state is asymptotically stable in the case $R_0<1$ and unstable
if $R_0>1$. In fact looking at the proof reveals that under the
assumptions of that theorem the following stronger statements hold. In the
case $R_0<1$ all eigenvalues of the linearization at the disease-free
steady state have negative real parts and in the case $R_0>1$ the
linearization has an eigenvalue with positive real part. This gives an
indirect proof that the inequalities
$\mathcal{A}_1 \mathcal{A}_2> \mathcal{A}_3$ and $\mathcal{A}_3>0$ of the last
section hold.

In order to understand the effects of vaccination it is useful to write
the basic reproductive ratio schematically in the form
$R_0=A\left(\frac{B\sigma+C}{D\sigma +E}\right)$ where
\begin{eqnarray}
&&A=\frac{\Lambda\xi}{\omega(\xi +\mu)(\mu+\delta+d)},\\
&&B=\epsilon \omega+\frac{\delta_{1n}\alpha_2 \varphi}{\kappa},\\
&&C=(\mu+t')\left(\omega\beta+\frac{\delta_{1n}\alpha_1 \varphi}{\kappa}\right),\\
&&D=\mu+\lambda t',\\
&&E=\mu(\mu+t').
\end{eqnarray}
The sign of the derivative of $R_0$ with respect to $\sigma$ is equal to that
of $BE-CD$. This last quantity is equal to
\begin{equation}
\omega (\mu+t')[(\epsilon-\beta)\mu-\lambda t'\beta]
 +\frac{\delta_{1n}\varphi}{\kappa}(\mu+t')[(\alpha_2-\alpha_1)\mu-\lambda t')].
\end{equation}
Under the assumptions made on the parameters it is negative and so we see
that increasing the vaccination rate decreases $R_0$, generalizing a result
of \cite{van2002reproduction}.

\textbf{3.2. Backward bifurcation analysis}

The concept of a backward bifurcation is used in the literature on
epidemiological models. It is defined in situations where a
definition of the basic reproduction ratio $R_0$ is available. In many models
endemic steady states only exist in the case $R_0>1$. We think of the direction
of increasing $R_0$ as the forward direction and that is where endemic steady
states occur. There are, however, models where it happens that near $R_0=1$
there are endemic steady states with $R_0<1$, i.e. in the backward direction.
The steady state bifurcates from the disease-free steady state as the parameter
$R_0$ is varied. A commonly occurring case is where there is a generic
transcritical bifurcation for $R_0=1$ and this covers both the forward and
backward cases, these being distinguished by the sign of a parameter $a$. We
call this case a generic backward bifurcation.
In our model the qualitative behaviour depends on a parameter which is a natural
number $n$. For $n\ne 2$ we show that any generic transcritical bifurcation
must be a forward bifurcation. Thus a generic backward bifurcation is
impossible in that case. For $n=2$ we show that generic backward bifurcations
do occur for some values of the parameters.

\noindent Consider the disease-free steady state $E^*=(S^*,0,0,V^*,R^*,0)$ for the
system \eqref{Eq:ModelExp1}-\eqref{Eq:ModelExp6}. We choose $\beta$ as the
bifurcation parameter and denote its value at the bifurcation point where
$R_0=1$ by $\beta^*$. Then
\begin{align}
    \beta^* =\frac{(\xi+\mu)(\mu+\delta+d)}{\xi S^*} - \frac{\delta_{1n}\alpha_1 \varphi}{\omega \kappa} -\left(\epsilon+\frac{\delta_{1n}\alpha_2 \varphi}{\omega \kappa}\right)\frac{\sigma}{t'+\mu},\label{Eq:bs} 
\end{align}
where $S^* =\dfrac{\Lambda (t'+\mu)}{\mu (\sigma+t'+\mu)+\lambda t' \sigma}.$
Note that for fixed values of the other parameters there is only a choice of
$\beta$ for which $R_0=1$ if the right hand side of (\ref{Eq:bs}) is
positive.

\noindent Following the ideas presented by \cite{van2002reproduction}, according to center manifold theory, it is necessary to compute right and left eigenvectors of the Jacobian matrix evaluated at the disease-free steady state $E^*$ and
$\beta=\beta^*$. Consider a right eigenvector of the form $\underline{w}=(w_1,w_2,w_3,w_4,w_5,w_6)^T$. Thus the system leads to

\begin{align}
  -(\sigma+\mu) w_1 - \beta^* S^* w_3 +(1-\lambda)t' w_4
  -\frac{\alpha_1\delta_{1n}}{\kappa} S^* w_6 &=0, \label{Eq:RE1}\\
  -(\xi+\mu) w_2 +(\beta^* S^*+\epsilon V^*) w_3
  + \frac{\delta_{1n}}{\kappa} (\alpha_1 S^*+\alpha_2 V^*)w_6 &=0,\label{Eq:RE2}\\
    \xi w_2 - (\delta+d+\mu) w_3&=0, \label{Eq:RE3}\\
  \sigma w_1 - \epsilon V^* w_3-(t'+\mu)w_4
  -\frac{\delta_{1n}\alpha_2}{\kappa}V^* w_6 &=0,\label{Eq:RE4}\\
    \delta \omega_3+\lambda t' w_4-\mu w_5 &=0,\label{Eq:RE5}\\
    \varphi w_3 -\omega w_6 &=0.\label{Eq:RE6}
\end{align}
Using Eq. \eqref{Eq:RE3} and \eqref{Eq:RE6}, we obtain
\begin{align}
    w_3=\dfrac{\xi w_2}{\delta+d+\mu} \quad \textrm{and} \quad w_6=\dfrac{\varphi \xi}{\omega(\delta+d+\mu)} w_2. \label{Eq:w3w6}
\end{align}
 Using these  in Eq. \eqref{Eq:RE2}, we find
\begin{align}
\left[    \beta^*- \frac{(\delta+d+\mu)(\xi+\mu)}{\xi S^*} + \frac{\delta_{1n}\alpha_1 \varphi}{\omega \kappa} +\left(\epsilon+\frac{\delta_{1n}\alpha_2 \varphi}{\omega \kappa}\right)\frac{\sigma}{t'+\mu}\right]w_2=0,
\end{align}
\noindent Thus we see that \eqref{Eq:bs} is a necessary condition for there to be a
vector in the kernel with $w_2\ne 0$. Note that if $w_2>0$ then $w_3$ and
$w_6$ are positive, as they must be as a consequence of the general theory.
Besides, using Eqs. \eqref{Eq:RE1} and \eqref{Eq:RE4}, one obtains
\begin{align}
    -(\sigma+\mu) w_1 +(1-\lambda)t' w_4 & = \left(\beta^* +\frac{\delta_{1n}\alpha_1 \varphi }{\omega \kappa}\right)\frac{\xi S^*}{(\delta+d+\mu)} w_2,\\
    \sigma w_1 - (t'+\mu) w_4 &= \left(\epsilon +\frac{\delta_{1n}\alpha_2 \varphi}{\omega \kappa}\right)\frac{\xi V^*}{(\delta+d+\mu)} w_2,
\end{align}
respectively. That leads to
\begin{align}
    w_1 = - \dfrac{\dfrac{\xi w_2 }{\delta+d+\mu}\left[\left(\beta^*+\dfrac{\delta_{1n}\alpha_1 \varphi}{\omega \kappa}\right)(t'+\mu)S^* + \left(\epsilon +\dfrac{\delta_{1n}\alpha_2 \varphi}{\omega \kappa}\right)(1-\lambda)t' V^*\right]}{\mu(\sigma+t'+\mu)+\sigma \lambda t'}<0, \label{Eq:w1}
\end{align}
for $\lambda \leq 1$ and
\begin{align}
    w_4 = - \dfrac{\dfrac{\xi w_2 }{\delta+d+\mu}\left[\left(\beta^*+\dfrac{\delta_{1n}\alpha_1 \varphi}{\omega \kappa}\right)\sigma S^* + \left(\epsilon +\dfrac{\delta_{1n}\alpha_2 \varphi}{\omega \kappa}\right)(\sigma+\mu) V^*\right]}{\mu(\sigma+t'+\mu)+\sigma \lambda t'}<0. \label{Eq:w4}
\end{align}
Furthermore using Eq. \eqref{Eq:RE5}:
\begin{align}
    w_5= \dfrac{\xi w_2 }{\delta+d+\mu} \left[ \delta-\lambda t' \dfrac{\left(\beta^*+\dfrac{\delta_{1n}\alpha_1 \varphi}{\omega \kappa}\right)\sigma S^* + \left(\epsilon +\dfrac{\delta_{1n}\alpha_2 \varphi}{\omega \kappa}\right)(\sigma+\mu) V^*}{\mu(\sigma+t'+\mu)+\sigma \lambda t'} \right].
\end{align}
In a similar manner, a left eigenvector can be written in the form $\underline{v} = (v_1,v_2,v_3,v_4,v_5,v_6)$ for which
\begin{align}
    v_5&=0, \label{Eq:LE1}\\
    -(\sigma+\mu)v_1+\sigma v_4 &=0,\label{Eq:LE2}\\
    -(\xi+\mu)v_2+\xi v_3 &=0,\label{Eq:LE3}\\
    -\beta^*S^* v_1 +(\beta^* S^*+\epsilon V^*) v_2-(\delta+d+\mu)v_3-\epsilon V^* v_4+\varphi v_6 &=0, \label{Eq:LE4}\\
    (1-\lambda)t' v_1 - (t'+\mu) v_4 &=0, \label{Eq:LE5}\\
    -\frac{\delta_{1n}\alpha_1}{\kappa} S^* v_1 +\frac{1}{\kappa} (\delta_{1n}\alpha_1 S^*+\delta_{1n}\alpha_2 V^*) v_2 - \frac{\delta_{1n}\alpha_2}{\kappa} V^* v_4 -\omega v_6 &=0. \label{Eq:LE6}
\end{align}
Using \eqref{Eq:LE2} and \eqref{Eq:LE5}, we find
\begin{align}
    \left(\sigma \lambda t' + \mu(\sigma +\mu+t')\right) v_4=0,
\end{align}
leading to $v_4=0$ and thus $v_1=0$. In addition, from Eqs. \eqref{Eq:LE3} and \eqref{Eq:LE6}, we obtain
\begin{align}
    v_3=\left(1+\frac{\mu}{\xi}\right) v_2,\quad \textrm{and} \quad v_6=\frac{\delta_{1n}}{\kappa \omega} \left(\alpha_1 S^* +\alpha_2 V^*\right) v_2.
\end{align}
Thus the left eigenvector becomes
\begin{align}
    \underline{v}=\left(0,v_2,\left(1+\frac{\mu}{\xi}\right)v_2,0,0,\frac{\delta_{1n}}{\kappa \omega}(\alpha_1 S^*+\alpha_2 V^*) v_2\right),\label{Eq:LET}
\end{align}
and we find
\begin{align}
    \underline{w} \cdot \underline{v} = \left(1+\frac{\xi+\mu}{\delta+d+\mu}+\frac{\delta_{1n}\varphi \xi}{\kappa \omega^2 (\delta+d+\mu)}\right) w_2 v_2 >0.
\end{align}
Let $a$ be the bifurcation coefficient introduced in \cite{van2002reproduction}. Considering the model \eqref{Eq:ModelExp1}-\eqref{Eq:ModelExp6} in the form $\dot{ x_i}=f_i(x_i),\;\; i=\{1,2,3,4,5,6\}$, it is given by
\begin{align}
    a= \sum\limits_{k,i,j=1}^{6} v_k w_i w_j \frac{\partial^2 f_k}{\partial x_i \partial x_j}(0,0). \label{Eq:Coeff}
\end{align}
Using Eqs. \eqref{Eq:LET} and \eqref{Eq:Coeff}
\begin{align}
    a =&v_2 \sum\limits_{i,j=1}^{6} w_i w_j \frac{\partial^2 f_2}{\partial x_i \partial x_j} + \left(1+\frac{\mu}{\xi}\right)v_2 \sum\limits_{i,j=1}^{6} w_i w_j \frac{\partial^2 f_3}{\partial x_i \partial x_j}\nonumber \\& + \frac{\delta_{1n}(\alpha_1 x_1^*+\alpha_2 x_4^*)}{\kappa \omega} v_2 \sum\limits_{i,j=1}^{6} w_i w_j \frac{\partial^2 f_6}{\partial x_i \partial x_j},\nonumber
\end{align}
where
\begin{align}
    f_2 &= \beta^* x_1 x_3 + \epsilon x_3 x_4 - \xi x_2-\mu x_2 + \alpha_1 x_1 g(x_6,\kappa)+\alpha_2 x_4g(x_6,\kappa) , \nonumber\\
    f_3 & = \xi x_2-(\delta+d+\mu)x_3,\nonumber \\
    f_6 & = \varphi x_3-\omega x_6.\nonumber
\end{align}
\noindent Note that $\frac{\partial g}{\partial C}(0,\kappa)=\kappa^{-1}$ for $n=1$ and
zero otherwise. $\frac{\partial^2 g}{\partial C^2}(0,\kappa)$ is equal to
$-2\kappa^{-2}$ for $n=1$, $2\kappa^{-1}$ for $n=2$ and zero otherwise.
Since second derivatives of $f_3$ and $f_6$ with respect to $x_i,\;\; i=\{1,2,3,4,5,6\}$ are always zero;
\begin{eqnarray}
  &&a=2 v_2 \left[w_1 w_3 \beta^* + w_3 w_4 \epsilon\right. \notag\\
  &&\left.+(\alpha_1w_1+\alpha_2w_4)w_6\frac{\partial g}{\partial C}(0,\kappa)
  +(\alpha_1x_1^*+\alpha_2x_4^*)
  w_6^2\frac{\partial^2 g}{\partial C^2}(0,\kappa)
  \right].
\end{eqnarray}
\noindent We want to determine the sign of $a$ and since $v_2>0$ this is the same
as that of the
expression in square brackets. The first two summands are negative. Now
$(w_1+w_4)w_6<0$ and $w_6^2>0$. Hence in the case $n=1$ we see that all
summands are negative and hence $a<0$. It follows that in that case the
conditions for a backward bifurcation given in \cite{van2002reproduction}
cannot be satisfied. The same conclusion is obtained in the case $n\ge 3$.
In the exceptional case $n=2$ the last summand is positive and so we
investigate further whether a backward bifurcation can take place in that case.
In the notation of \cite{van2002reproduction} we have
\begin{align}
    b=   \sum\limits_{k,i=1}^{6} w_k w_i \frac{\partial^2 f_k}{\partial x_i \partial \beta}(0,0). \label{Eq:Coeff2}
\end{align}
Since only $f_1$ and $f_2$ involve the parameter $\beta$, derivatives with
respect to $x_1$ and $x_3$ are non-zero. Here
\begin{align}
    b=w_3 x_1^* (w_2-w_1) >0.
\end{align}
It follows from \cite{van2002reproduction} that there is a backward bifurcation
precisely when the right hand side of (\ref{Eq:bs}) is positive and
\begin{equation}\label{bthreshold}
2(\alpha_1S^*+\alpha_2V^*)w_6^2\kappa^{-1}>-w_3(w_1\beta^*+w_4\epsilon).
\end{equation}
Note that in the case $n>1$ the expression for $\beta^*$ simplifies to
\begin{align}\label{betastar}
    \beta^* =\frac{(\xi+\mu)(\mu+\delta+d)\left[\mu(\sigma+t'+\mu)+\lambda t' \sigma\right]-\Lambda \xi \epsilon \sigma}{\Lambda \xi (t'+\mu)}.
\end{align}
Do there exist values of the parameters for which these conditions are
satisfied? To investigate this we substitute the expressions for $w_6$, $w_3$,
$w_1$ and $w_4$ into (\ref{bthreshold}). The result is
\begin{eqnarray}
&&  2(\alpha_1S^*+\alpha_2V^*)
   \dfrac{\varphi^2}{\omega^2\kappa}
   \nonumber\\
  &&  >\frac{\Lambda[(\beta^*)^2(t'+\mu)^2+
     \beta^*\epsilon((2-\lambda)t'+\mu)\sigma
     +\epsilon^2t'\sigma (\sigma+\mu)]}
     {[\mu(\sigma+t'+\mu)+\sigma\lambda t']^2}.
\end{eqnarray}
Using the expressions for $S^*$ and $V^*$ this can be simplified to
\begin{eqnarray}
&&  2(\alpha_1(t'+\mu)+\alpha_2\sigma)
   \dfrac{\varphi^2}{\omega^2\kappa}
   \nonumber\\
  &&  >\frac{[(\beta^*)^2(t'+\mu)^2+
     \beta^*\epsilon((2-\lambda)t'+\mu)\sigma
     +\epsilon^2t'\sigma (\sigma+\mu)]}
     {\mu(\sigma+t'+\mu)+\sigma\lambda t'}.
     \label{Eq:backwardbif}
\end{eqnarray}
\noindent Thus it is clear that if $\alpha_1$ or $\alpha_2$ is made large enough while the other parameters are kept fixed then there is a backward bifurcation. It is important to note that the parameters given for the case $n=2$ in Table \ref{tab:my_label2} satisfy the conditions for a backward bifurcation in \eqref{Eq:backwardbif}. These results with now be summed up.

\noindent
{\bf Theorem 3.1} If $n=2$ and the parameters in the system
\eqref{Eq:ModelExp1}-\eqref{Eq:ModelExp6} satisfy the inequality
\ref{Eq:backwardbif} with the quantity $\beta^*$ defined by \ref{betastar}
being positive then the parameter $a$ of \cite{van2002reproduction} is
positive and a generic backward bifurcation occurs. There exist parameters
for which these conditions are satisfied. If $n\ne 2$ the condition $a>0$ is
never satisfied.


\noindent The centre manifold at the bifurcation point is one-dimensional. Since $v_3>0$ we can use $I$ as a parameter on the centre manifold. The restriction of the dynamical system to the centre manifold is of the form $\dot I=f(I,\beta)$, where we have suppressed the dependence on the parameters other than $\beta$ in the notation. With this notation the sign of $a$ is equal to that of $\frac{\partial^2f}{\partial I^2}$ while that of $b$ is equal to that of $\frac{\partial^2f}{\partial I\partial\beta}$. We have $f(0,\beta)=0$ for all $\beta$. Moreover there is a curve $c(\beta)$ of steady states with
$c(\beta^*)=0$. The sign of $c'(\beta^*)$ is equal to that of $a$.

\textbf{4. Endemic steady states}

\noindent In this section we consider endemic steady states, i.e. those where all
the unknowns are positive. It follows from (\ref{Eq:ModelExp4}) that at a
steady state $V=\frac{\sigma}{\epsilon I+(t'+\mu)+\alpha_2 g}S$. Substituting
this into (\ref{Eq:ModelExp1}) gives
\begin{eqnarray}\label{SofI}
&&  S=\frac{\Lambda}{\beta I+\sigma+\alpha_1g+\mu}
  \left[1-\frac{(1-\lambda)t'\sigma}
   {(\beta I+\sigma+\alpha_1 g+\mu)(\epsilon I+(t'+\mu)+\alpha_2 g)}\right]^{-1}
   \nonumber\\
  &&  =\frac{\Lambda(\epsilon I+(t'+\mu)+\alpha_2 g(C))}
     {(\beta I+\sigma+\alpha_1 g(C)+\mu)(\epsilon I+(t'+\mu)+\alpha_2 g(C))
     -(1-\lambda)t'\sigma}
\end{eqnarray}
\noindent Note also that due to (\ref{Eq:ModelExp6}) we have $C=\frac{\varphi}{\omega}I$
and that due to (\ref{Eq:ModelExp3}) we have $E=\frac{\xi}{\delta+d+\mu}I$.
These relations allow $S$, $E$, $V$ and $C$ to be expressed in terms of $I$.
Thus substituting them into (\ref{Eq:ModelExp2}) gives an equation for $I$
alone. Each summand contains a factor $I$ and for an endemic steady state this
can be cancelled. There remains
\begin{equation}
0=\beta S+\epsilon V-\frac{(\xi+\mu)(\delta+d+\mu)}{\xi}
     +(\alpha_1S+\alpha_2V)\frac{(\varphi/\omega)^nI^{n-1}}{(\varphi/\omega)^nI^n
     +\kappa}\label{prepol}
\end{equation}
\noindent Denote the expression in the denominator of (\ref{SofI}) by $Z$.
Multiplying (\ref{prepol})
by $Z(\epsilon I+(t'+\mu)+\alpha_2 g)[(\varphi/\omega)^nI^n+\kappa]$
gives a polynomial equation $p(I)=0$ for $I$. Endemic steady states are in one
to one correspondence with positive roots of $p$.

\noindent Since we are most interested in the case with backward bifurcations
we now restrict to the case $n=2$, where there are considerable simplifications
in these formulae. It is clear that $p(I)\to -\infty$ as $I\to\infty$. Moreover
$p(0)$ is equal to a positive factor times
\begin{equation}
\beta S^*+\epsilon V^*-\frac{(\xi+\mu)(\delta+d+\mu)}{\xi}.
\end{equation}
\noindent We see that the sign of $p(0)$ is the same as that of $R_0-1$. As $\beta$
increases through $\beta^*$ the sign of $p(0)$ changes from negative to
positive. When $p'(0)\ne 0$ we have a generic transcritical bifurcation
(cf. \cite{wiggins90}, Section 3.1).
The sign of $p'(0)$ is the same as that of $c'(\beta^*)$ in the
discussion of the centre manifold in the previous section.
If $p'(0)<0$ then for $\beta$ slightly greater than $\beta^*$ there
exists a positive root of $p$ close to zero and hence a positive steady
state close to the disease free steady state. This corresponds to a forward bifurcation. If, on the other hand, $p'(0)>0$ then for $\beta$ slightly less than $\beta^*$ there exists a positive root of $p$ close to zero and hence a
positive steady state close to the disease free steady state. This corresponds
to a backward bifurcation. In this case $p$ is positive for $I$ slightly
larger than its value $I_1$ at that steady state. By the intermediate value
theorem there must exist some $I_2>I_1$ with $p(I_2)=0$ and hence a second
positive steady state. The direction of the flow on the centre manifold shows
that the steady state with $I=I_1$ is unstable. The stability of the steady
state with $I=I_2$ cannot be determined by the arguments we have presented.

\textbf{5. Simulations}

\noindent In a broad context, the process of mathematical modelling and data fitting revolves around formulating mathematical models that describe real-world phenomena and then adjusting the parameters of these models to best match observed data. Therefore, the aim is to capture the underlying relationships and behavior of the system being studied and use the available data to validate the model. In our study, we ensure the validation of the mathematical model of the COVID-19 outbreak by using the data fitting of the model regarding the observed data. However, the available data is scarce (only the actual data for infectious and vaccinated people are almost certain). Vaccination is globally considered to be the most effective solution for infectious diseases such as the recent COVID-19 outbreak. In this section, as an example, model results and observed data for the vaccinated class are compared for COVID-19 scenarios in Turkey.  Some of the realistic parameters are taken from the literature, see the related references in Table \ref{tab:my_label2} for a detailed discussion of parameter choices. The remaining seven parameters are estimated, by fitting the vaccinated compartment generated from the system \eqref{Eq:ModelExp1}-\eqref{Eq:ModelExp6} to the observed number of COVID-19 vaccinated individuals using standard model-fitting procedures.

\noindent The least squares method is the process of finding the best-fitting curve or line of best fit for a set of data points by reducing the sum of the squares of the offsets (residual part) of the points from the curve.
The vector consisting of seven parameters $p=(\beta,\epsilon,\alpha_1,\alpha_2,\xi,\delta,\sigma)$ can be estimated via \textit{parameter estimation}. In this context, the model given by \eqref{Eq:ModelExp1}-\eqref{Eq:ModelExp6}  is evaluated by considering a non-linear least squares problem with positive constraints, where the best fitting curve can be found for a small data set of vaccinated class by minimising the sum of squares of the deviations of data points from the plotted graph \cite{bulut2021modelling,coleman1996interior}. This may be described as
\begin{equation}
    S(p) = \sum \left(V_i-\mathcal{F}(x_i,p)\right)^2,
\end{equation}
where $V_i$ represents the data set and $\mathcal{F}(x_i,p)$ denotes the model result with for a vector of unknowns $p$. To minimise the function $S(p)$, the non-linear least square minimization routine \textit{lsqcurvefit} of MATLAB is used \cite{MATLAB}. Parameters obtained from this approach are given in Table \ref{tab:my_label2}. Besides MATLAB’s standard \textit{ode45} solver \cite{MATLAB} is applied for numerical integration of the system \eqref{Eq:ModelExp1}-\eqref{Eq:ModelExp6} with suitable initial conditions provided in the text.

\begin{figure}{h}
    \centering
    \includegraphics[scale=0.5]{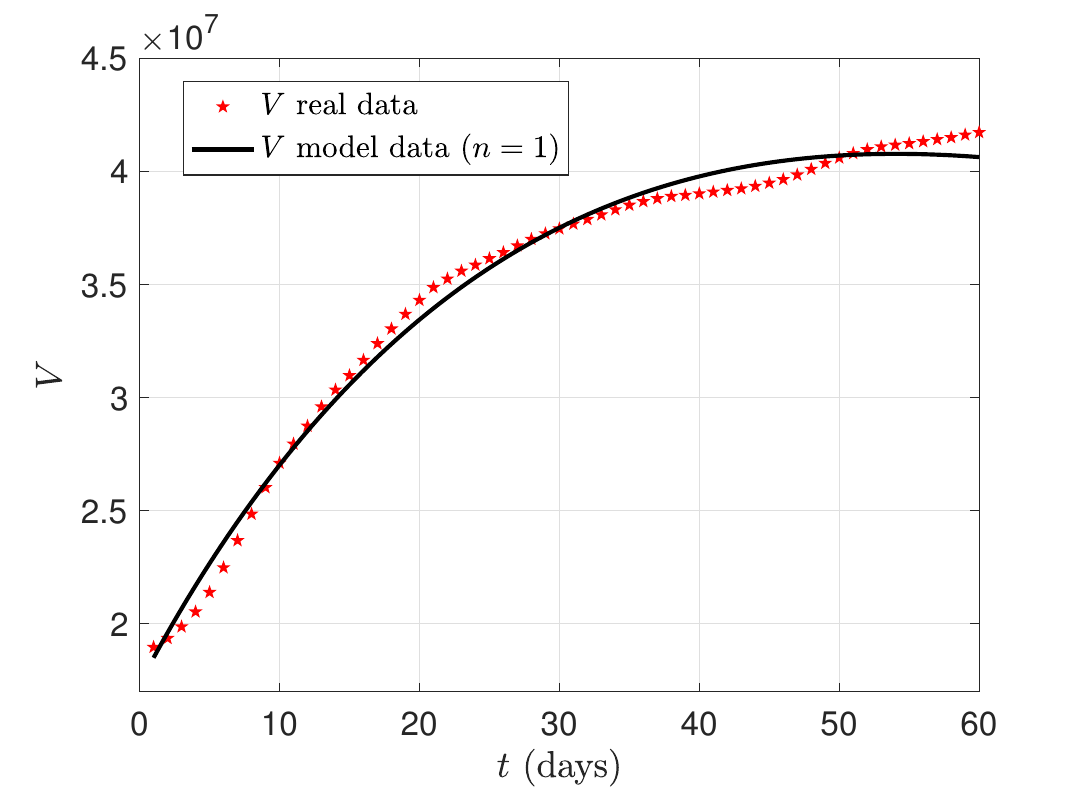}
    \caption{Vaccination component for the model for the case with $n=1$ is compared with real data for the period of June $10$, $2021$ and August $8$, $2021$. Parameters are given in Table \ref{tab:my_label2}. Initial conditions are chosen as $S_0=61098000, V_0=18500000, E_0=2200000, R_0=2000, C_0=20000$. The black line denotes the model result and the red stars represent daily vaccinated cases.  }
    \label{fig:S2V}
\end{figure}

\noindent Numerical simulations of the model \eqref{Eq:ModelExp1}-\eqref{Eq:ModelExp6} can be shown with the parameters given in Table \ref{tab:my_label2} for $n=1$. An example data set of vaccinated people during the COVID$-19$ outbreak in Turkey is taken from the World Health Organisation. In Fig. \ref{fig:S2V}, the results are shown for model \eqref{Eq:ModelExp1}-\eqref{Eq:ModelExp6} with $n=1$  fitted to the data of individuals vaccinated between June $10$, $2021$ and August $8$, $2021$. The total population is taken $83$ million \cite{TUIK}.

\begin{figure}
    \centering
    \includegraphics[scale=0.5]{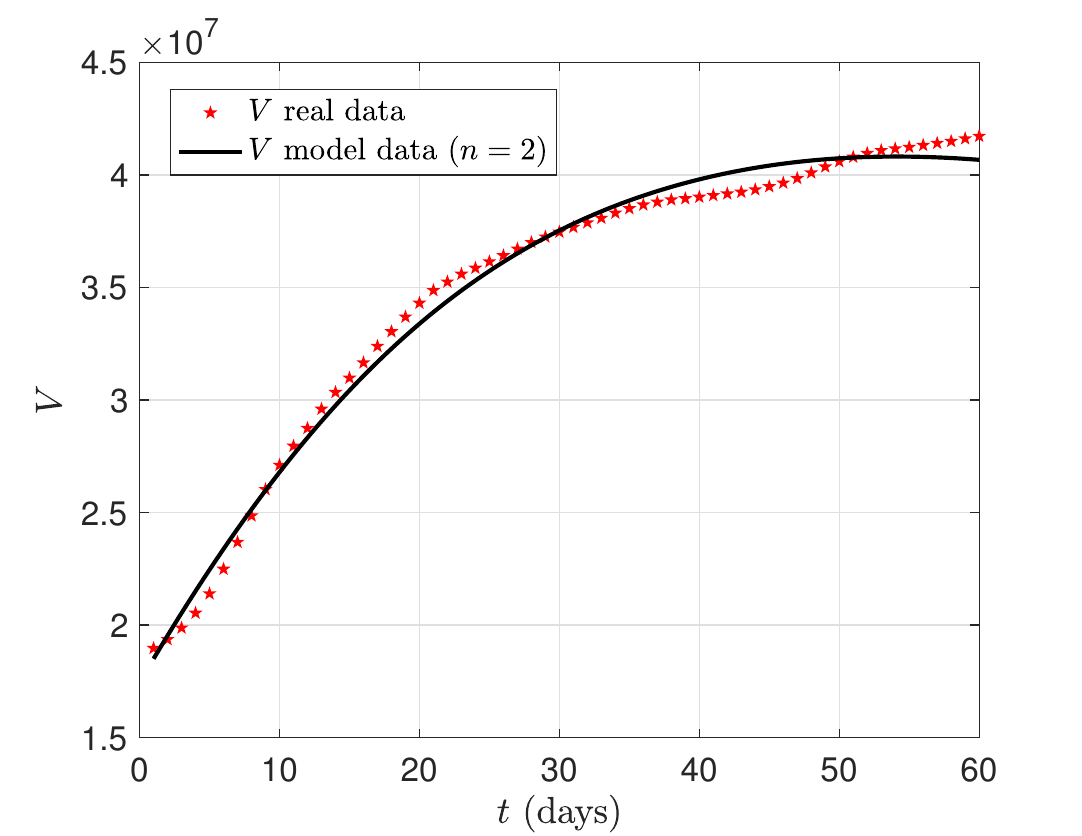}
    \caption{Vaccination component for the model for the case with $n=2$ is compared with real data for the period of June $10$, $2021$ and August $8$, $2021$. Parameters are given in Table \ref{tab:my_label2}. Initial conditions are chosen as $S_0=61098000, V_0=18500000, E_0=2200000, R_0=2000, C_0=20000$. The black line denotes the model result and the red stars represent daily vaccinated cases.  }
    \label{fig:S4V}
\end{figure}

\noindent In Fig. \ref{fig:S4V}, success of parameter estimation is again demonstrated. Here numerical simulations of the model \eqref{Eq:ModelExp1}-\eqref{Eq:ModelExp6} with $n=2$ is shown with the parameters given in Table \ref{tab:my_label2} and the resulting outcome is compared with the data set of vaccinated people between June $10$, $2021$ and August $8$, $2021$. As seen, in Figs. \ref{fig:S2V} and \ref{fig:S4V},  the black curve corresponding to the model becomes flattened after around day $40$ and this agrees with the real data, where the real number of vaccinated people rapidly increases from $18967237$ to $41726338$ and it is smoothed roughly about July $20$.

\begin{table}[]
    \centering
       \caption{Estimated parameters of the model.}
    \begin{tabular}{c l l c l}
       \hline
       Parameters  &  Value ($n=1$) & Value ($n=2$) & Unit  & Source\\
       \hline
       $\Lambda$ & $3032$ & $3032$  & day$^{-1}$ & assumed based on  \cite{bajiya2020mathematical,garba2013cross,gumel2006sveir} \\
        $\beta$ & $1.0257\times 10^{-8}$ & $1.004 \times 10^{-8}$  & day$^{-1}$ & estimated  \\
        $\mu$ & $3.653 \times 10^{-5}$ & $3.653 \times 10^{-5}$  & day$^{-1}$ & assumed based on \cite{TUIK,bulut2021modelling,gumel2006sveir} \\
        $\epsilon$ & $1 \times 10^{-8}$ & $1 \times 10^{-8}$& day$^{-1}$ & estimated  \\
        $t'$ & $0.0055$  & $0.0055$  & day$^{-1}$ & assumed based on \cite{angeli2022modeling} \\
        $\lambda$ & $0.8$  & $0.8$ & day$^{-1}$ & assumed based on \cite{angeli2022modeling,baden2021efficacy,polack2020safety}  \\
        $d$     &  $0.1$   & $0.1$   & day$^{-1}$ & assumed based on \cite{bulut2021modelling,ritchie2020mortality} \\
        $\alpha_1$ &  $0.00041$ & $0.0001$ & day$^{-1}$ & estimated \\
        $\alpha_2$ & $0.00031$  &  $0.0001$ & day$^{-1}$ & estimated  \\
        $\omega$ & $5$ & $5$ & day$^{-1}$ &  \cite{bulut2021modelling}  \\
        $\varphi$ & $2$ & $2$ & day$^{-1}$ &  \cite{bulut2021modelling} \\
        $\kappa$ & $20000$ & $20000$ & copies $/$day & assumed based on \cite{drosten2004evaluation,bulut2021modelling} \\
        $\xi$ & $0.01004 $   & $0.01133$ & day$^{-1}$ & estimated \\
        $\delta$ & $0.19999$ & $0.2$& day$^{-1}$ &estimated \\
        $\sigma$ & $0.02136$  & $0.02126$ & day$^{-1}$ &estimated \\
    \end{tabular}
    \label{tab:my_label2}
\end{table}

\textbf{6. Sensitivity analysis}

\noindent Since varying the parameters may have a significant impact on the model output, one can perform a sensitivity analysis of the dynamical model to determine which model parameters are more influential for the dynamics. The parameters associated with the basic reproduction ratio $R_0$ have particular importance for the robustness of the model. In this context, the aim of the sensitivity analysis is to identify the most substantial parameter in the model for disease transmission.

\noindent Following the ideas presented in \cite{chitnis2008determining,martcheva2015introduction}, the sensitivity analysis  can be performed based on the basic reproduction ratio as
\begin{equation}
  \mathcal{S}_{i}^{\mathcal{P}}  = \frac{\mathcal{P}}{R_0} \frac{\partial R_0}{\partial \mathcal{P}},
\end{equation}
where $\mathcal{P}$ represents a generic parameter in the model \eqref{Eq:ModelExp1}-\eqref{Eq:ModelExp6}.

\begin{figure}[ht!]
\centering
\begin{tabular}{cccc}
\includegraphics[width=0.45\textwidth]{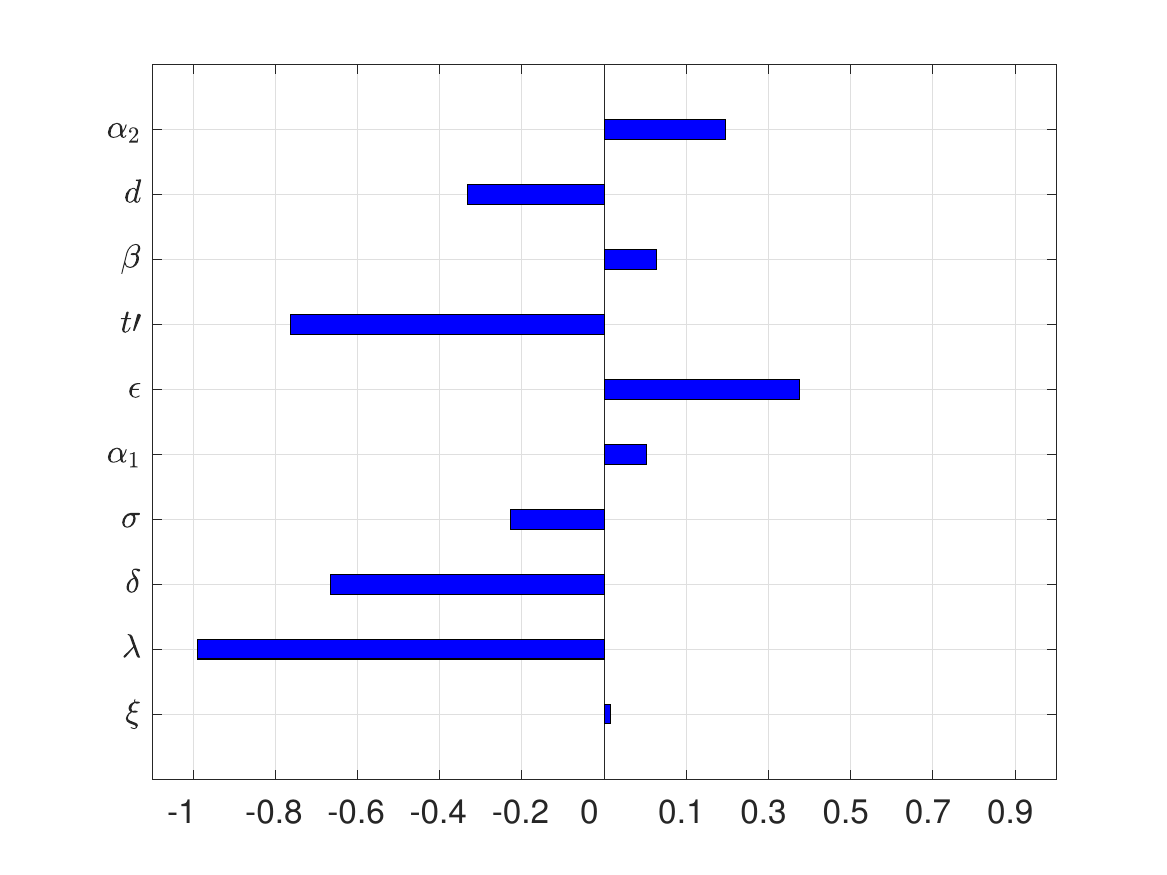} \\
\textbf{(a)}  \\ [6pt]
\end{tabular}
\hspace{-1cm}
\begin{tabular}{cccc}
\includegraphics[width=0.45\textwidth]{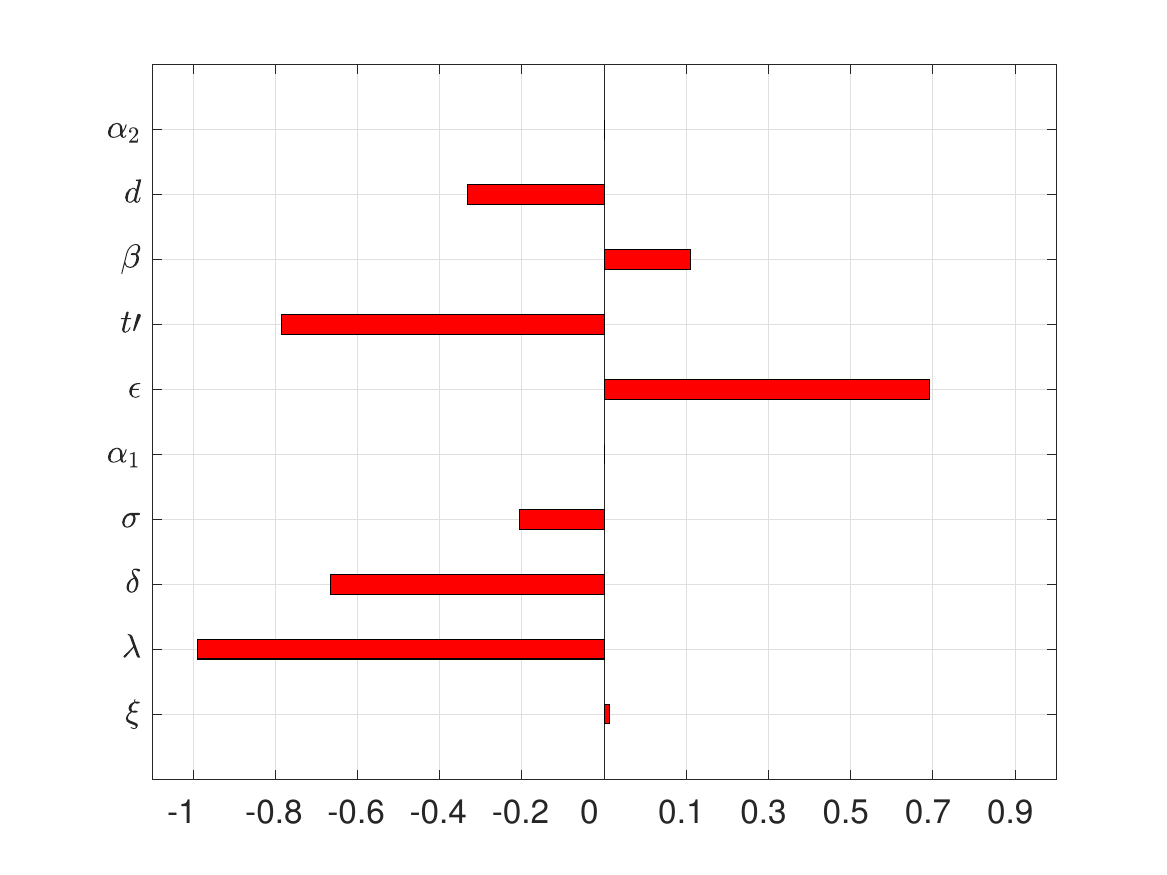} \\
\textbf{(b)}\\[6pt]
\end{tabular}
\caption{Sensitivity indices based on basic reproduction ratio $R_0$ with respect to various parameters given for Cases $n=1$ (a) and $n=2$ (b) in the  model \eqref{Eq:ModelExp1}-\eqref{Eq:ModelExp6} for Turkey. Parameter values given in Table \ref{tab:my_label2} are considered for both cases. }
\label{Fig:Sensitivity}
\end{figure}

\noindent The sensitivity indices of the system parameters in
Table \ref{tab:my_label2}  are demonstrated in \ref{tab:my_label3}, and also in Figs. \ref{Fig:Sensitivity}(a) for $n=1$ and \ref{Fig:Sensitivity}(b) for $n=2$.

\noindent As seen from Table \ref{tab:my_label3} and Fig. \ref{Fig:Sensitivity}, the model \eqref{Eq:ModelExp1}-\eqref{Eq:ModelExp6} is highly sensitive to the
parameters $\lambda$, $t'$ and $\delta$. Thus it can be concluded that an increase in these parameters  diminishes the basic reproduction ratio $R_0$ for both $n=1$ and $n=2$. The significance of some parameters may be different between
Case $1$ and Case $2$. For example, although the increase in the parameter $\epsilon$, i.e. the rate at which a vaccinated individual becomes exposed after being in contact with an infected individual, leads to an essential stimulus for the basic reproduction ratio for Case $2$, yet it has a much smaller impact on the $R_0$ for the case $n=1$.

\begin{table}[]
    \centering
       \caption{Sensitivity indices of basic reproduction ratio of model \eqref{Eq:ModelExp1}-\eqref{Eq:ModelExp6}, considering at the baseline parameters provided in Table \ref{tab:my_label2}.
number for $n=1$ and $n=2$}
    \begin{tabular}{c l l}
       \hline
       Parameters  &  Index ($n=1$) &  Index ($n=2$) \\
       \hline
        $\beta$ & $0.1276764$ & $0.2089141$   \\
        $\epsilon$ & $0.4754643$ & $0.7910858$ \\
        $t\prime$ &  $-0.7636042$ & $-0.7843018$  \\
        $\lambda$ & $-0.9897355$ & $-0.9897256$  \\
        $d$ &   $-0.3333038$    & $ -0.3332927$  \\
        $\alpha_1$ & $0.1020714$ & $0$\\
        $\alpha_2$ & $0.2947878$ & $0$ \\
        $\xi$ &   $0.0036252$  & $0.0032138$\\
        $\delta$ & $-0.6665743$  & $-0.6665855$\\
        $\sigma$ &  $-0.2276182$ & $-0.2067744$\\
    \end{tabular}
    \label{tab:my_label3}
\end{table}

\begin{figure}[ht!]
\centering
\begin{tabular}{cccc}
\includegraphics[width=0.45\textwidth]{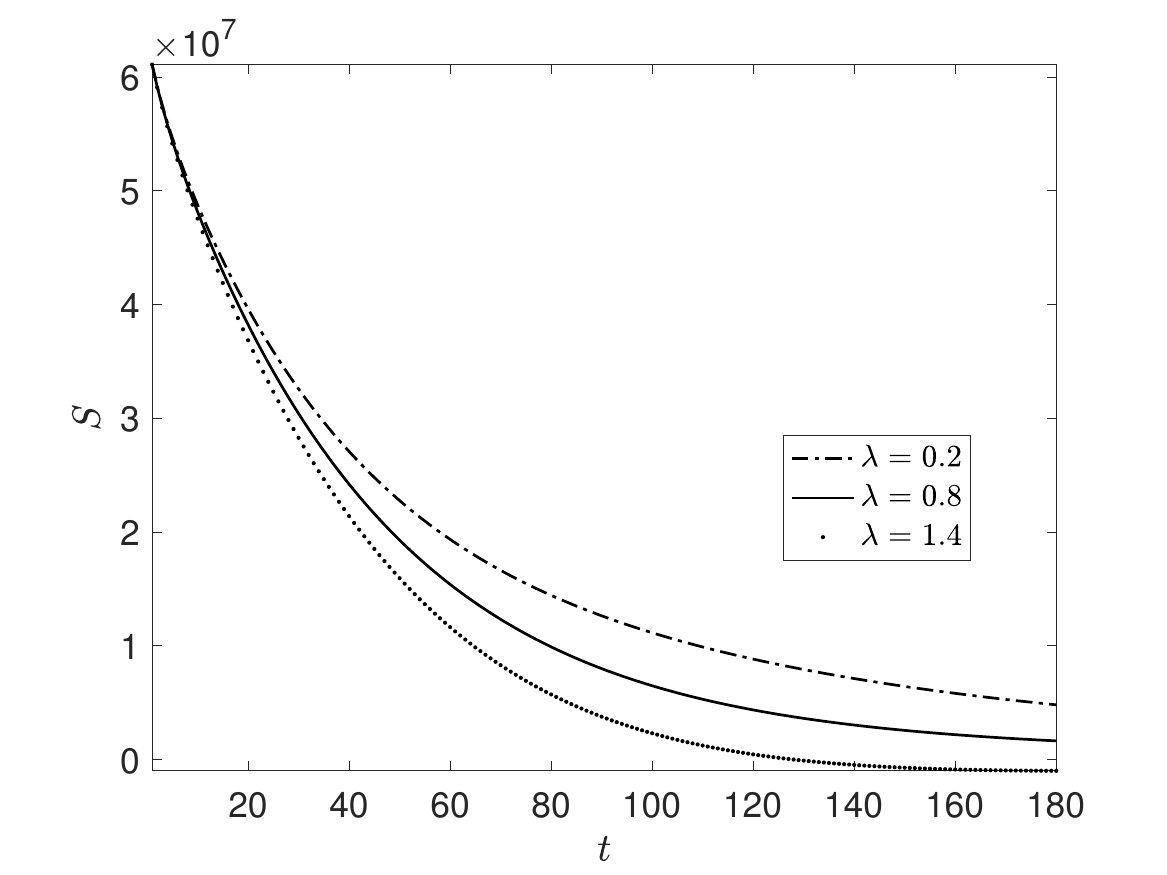} \\
\textbf{(a)}  \\ [6pt]
\end{tabular}
\hspace{-1cm}
\begin{tabular}{cccc}
\includegraphics[width=0.45\textwidth]{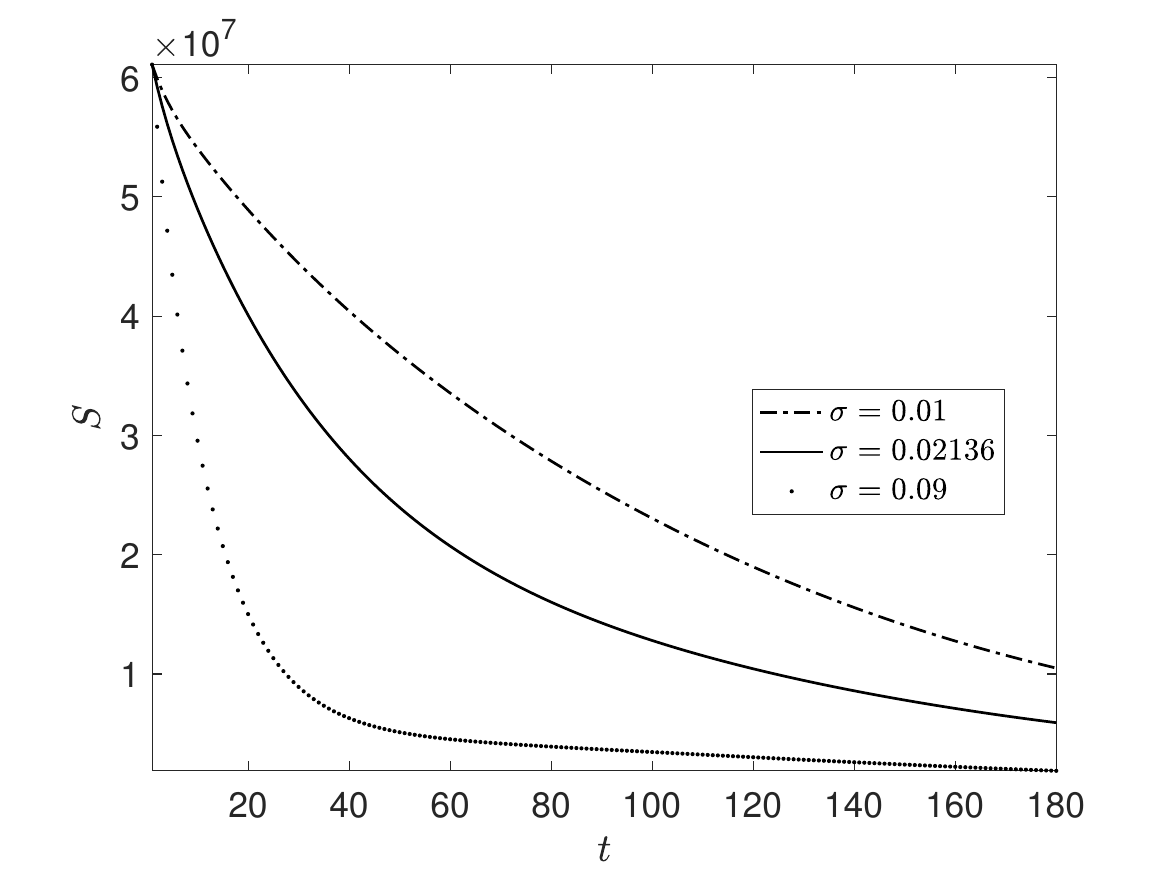} \\
\textbf{(b)}\\[6pt]
\end{tabular}
\caption{Plots for Susceptible population compartment for different  values of vaccine efficiency $\lambda=0.2,0.8,1.4$ and for different values of vaccination rate of susceptible individuals after first shot, $\sigma=0.01, 0.02136, 0.09$ with $n=1$. }
\label{Fig:S}
\end{figure}

\noindent Since the aim of our work is to further broaden the current knowledge of the modelling of recent  COVID-19 outbreak  with  vaccination, we focus on the role of two parameters regarding vaccination in the model.
In Fig. \ref{Fig:S}, time simulations of the Susceptible compartment over a period of $180$ days are presented for various rates of $\lambda=0.2,0.8,1,4$ (a) and $\sigma=0.01,0.02136,0.09$ (b), respectively denotes the efficiency of the vaccine and the vaccination rate of susceptible (after the first shot). As is seen from Fig. \ref{Fig:S} (a,b), with  the increase in $\lambda$ and $\sigma$,  susceptible individuals diminish at an earlier time and enter the vaccinated class.

\begin{figure}[ht!]
\centering
\begin{tabular}{cccc}
\includegraphics[width=0.45\textwidth]{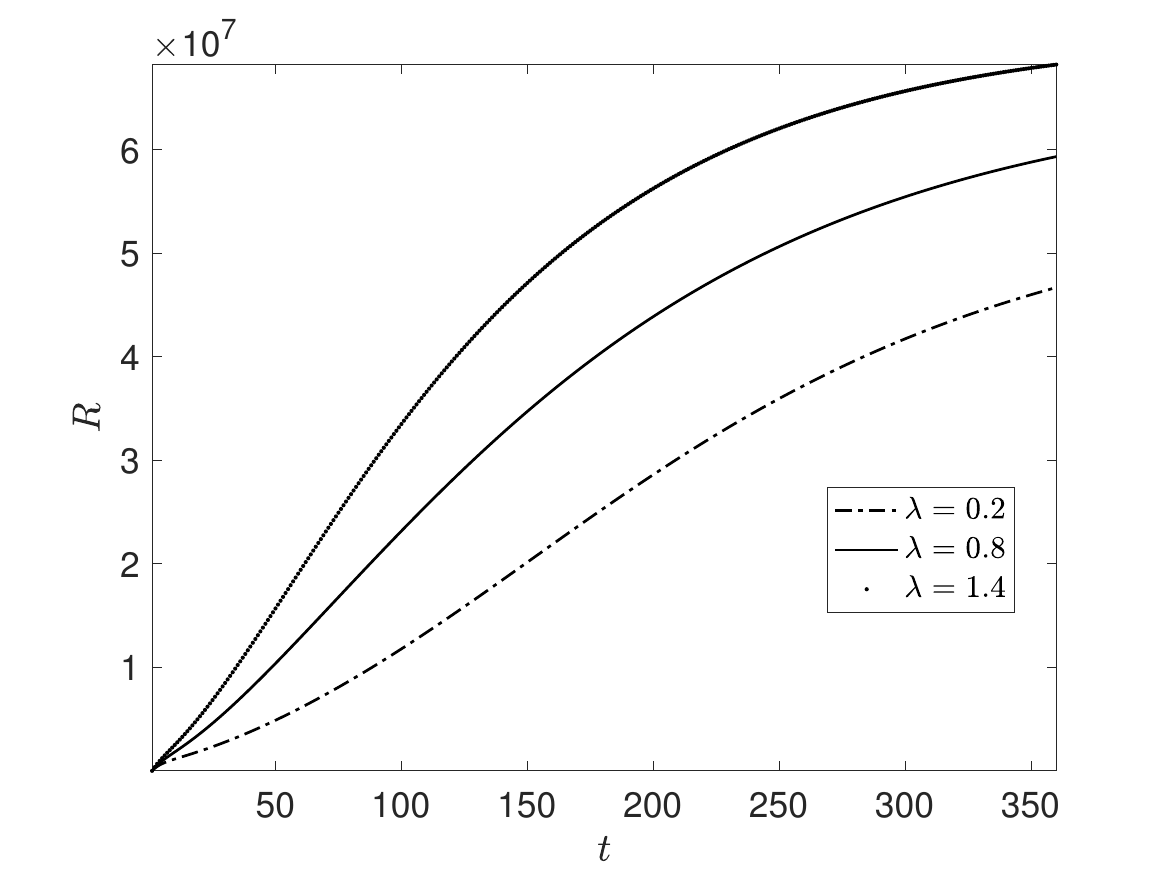} \\
\textbf{(a)}  \\ [6pt]
\end{tabular}
\hspace{-1cm}
\begin{tabular}{cccc}
\includegraphics[width=0.45\textwidth]{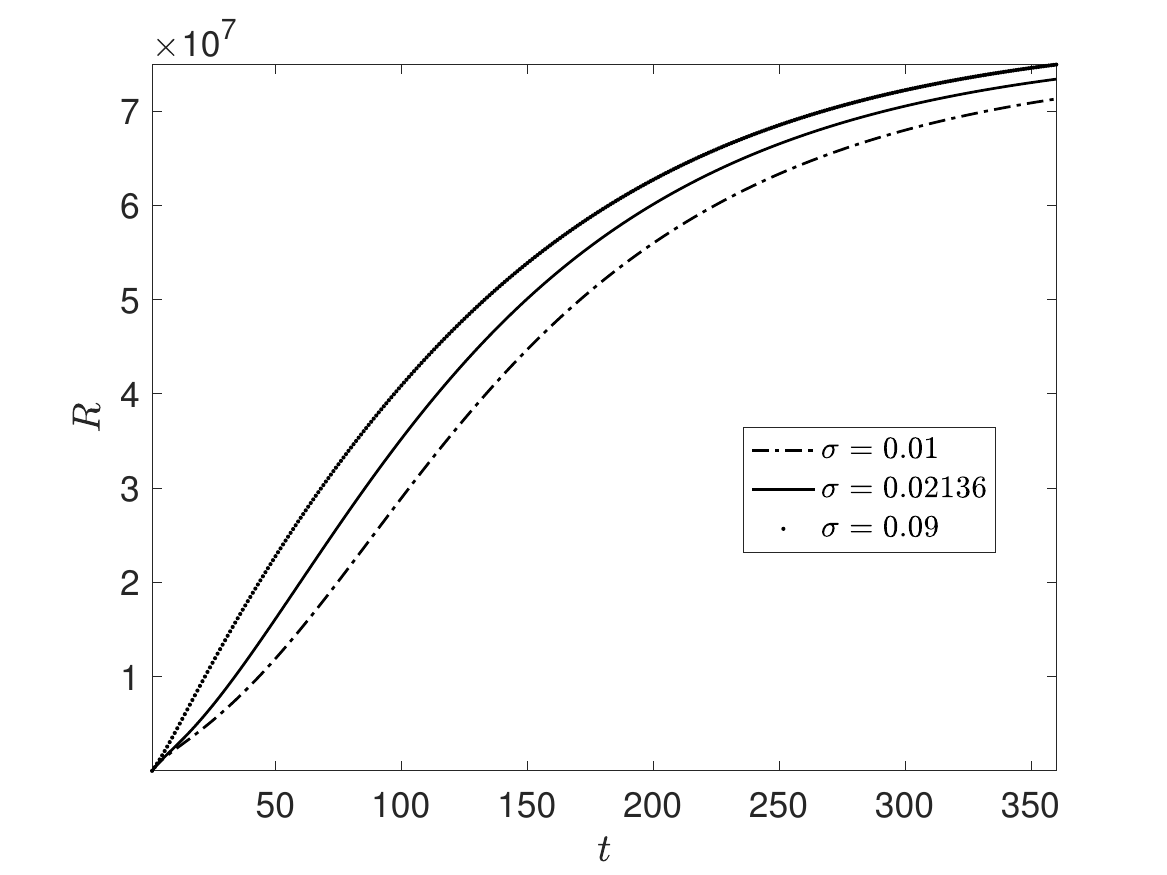} \\
\textbf{(b)}\\[6pt]
\end{tabular}
\caption{Plots for Recovery compartment for different  values of vaccine efficiency $\lambda=0.2,0.8,1.4$; and for different values of vaccination rate of susceptible individuals after first shot, $\sigma=0.01, 0.02136, 0.09$ with $n=1$.}
\label{Fig:R}
\end{figure}

\noindent Figure \ref{Fig:R} demonstrates time simulations of the recovered class over a period of $360$ days for various rates of $\lambda=0.2,0.8,1,4$ and $\sigma=0.01,0.02136,0.09$.

\textbf{Conclusion and outlook}

\noindent In this paper a model for an epidemic with a partially effective vaccination
and infection by virus in the environment is studied. Two different
implementations of the idea of a partially effective vaccination are included
and which of these is chosen does not seem to have an essential influence on
the qualitative behaviour of the solutions. In modelling infections coming
from the environment we used a phenomenological model for the force of infection
containing an integer $n\ge 1$ as a parameter. We discovered that the choice
$n=2$ leads to the occurrence of backward bifurcations while other choices of
$n$ do not. Thus here there is a major qualitative difference. On the other
hand, fitting to real data for COVID-19 in Turkey showed that both the cases
$n=1$ and $n=2$ worked and there was no clear indication that one of these
choices was better than the other according to that criterion.

\noindent In the case $n=2$ it was shown that for certain values of the parameters there
exists an endemic steady state although $R_0<1$. It was also shown that in
this situation there exist more than one endemic steady state. At least one of
the positive steady states is unstable. This confirms rigorously that certain
aspects of the usual picture of a backward bifurcation are present in this
model. An aspect of this picture which was not reproduced here is that the other
positive steady state should be stable. It would be desirable to prove a
stability statement of this kind analytically.

\noindent An outstanding challenge is to provide a mechanistic derivation of the response
function for infections coming from the environment. If this could be done
then it would help to decide which value of $n$ in the Ansatz we used is more
appropriate for modelling a given disease or whether, indeed, a different type
of function would give better results.

\let\thefootnote\relax
\footnotetext{\textbf{Funding:}
Ayt\"{u}l G\"{o}k\c{c}e was partially supported by a grant from the Niels Hendrik Abel Board.}

\footnotetext{{$^*$}(Corresponding author)}
\end{document}